\documentclass[proceedings,submission,pdftex]{dmtcs}
\usepackage[utf8x]{inputenc}
\usepackage[T1]{fontenc}
\usepackage{ae,aecompl}
\usepackage{graphics}

\usepackage{amssymb}
\usepackage{amsmath}

\usepackage{tikz}
\usetikzlibrary{matrix,shapes}

\usepackage{enumerate}
\usepackage{verbatim}
\usepackage{hyperref}
\usepackage{multicol}

\usepackage{ifthen}
\usepackage{verbatim}
\newboolean{draft}
\setboolean{draft}{false}
\newboolean{longversion}
\setboolean{longversion}{false}

\newtheorem{theorem}{Theorem}[section]
\newtheorem{lemma}[theorem]{Lemma}
\newtheorem{proposition}[theorem]{Proposition}
\newtheorem{corollary}[theorem]{Corollary}
\newtheorem{definition}[theorem]{Definition}
\newtheorem{example}[theorem]{Example}
\newtheorem{computerExample}[theorem]{Sage Example}

\newtheorem{conjecture}[theorem]{Conjecture}

\usepackage{xspace}
\newcommand{\contracting}{contracting\xspace}

\newcommand{\cuteq}{\sqsubseteq}
\newcommand{\fiber}{\operatorname{fiber}}

\newcommand{\End}{\operatorname{End}}

\newcommand{\one}{1}
\newcommand{\im}{\operatorname{im}}
\newcommand{\K}{\mathbb{K}}

\newcommand{\ls}[1][s] {{\overleftarrow{#1}}}
\newcommand{\N}{\mathbb{N}}
\newcommand{\opi}{\overline{\pi}}

\newcommand{\rank}{\operatorname{rank}}
\newcommand{\sg}[1][n]{{\mathfrak{S}_{#1}}}
\newcommand{\suchthat}{\mid}

\newcommand{\trans}{\operatorname{trans}}
\newcommand{\type}{\operatorname{type}}
\newcommand{\Z}{\mathbb{Z}}

\newcommand{\sym}{{\operatorname{Sym}}}
\newcommand{\qsym}{{\operatorname{QSym}}}
\newcommand{\fqsym}{{\operatorname{FQSym}}}

\newcommand{\lfix}{\operatorname{lfix}}
\newcommand{\rfix}{\operatorname{rfix}}

\newcommand{\lcoset}[2]{{}^{#2}\!#1}
\newcommand{\rcoset}[2]{\!#1^{#2}}

\newcommand{\biheckemonoid}[1][]{M}
\newcommand{\biheckealgebra}[1][W]{\mathcal{H}#1}
\newcommand{\wbiheckealgebra}[1]{\mathcal{H}W^{(#1)}}
\newcommand{\Ks}[1]{\mathcal{B_R}(#1)}
\newcommand{\RKs}[1]{\mathcal{RB_R}(#1)}
\newcommand{\Js}[1]{\mathcal{B_L}(#1)}
\newcommand{\RJs}[1]{\mathcal{RB_L}(#1)}

\newcommand{\Des}{{\operatorname{D}_R}}
\newcommand{\Rec}{{\operatorname{D}_L}}
\newcommand{\Jblock}[1]{J^{(#1)}}
\newcommand{\Kblock}[1]{K^{(#1)}}

\newcommand{\join}{\vee}

\ifdraft
\newcommand{\TODO}[2][To do: ]{\textcolor{red}{\textbf{#1#2}}}
\else
\newcommand{\TODO}[2][]{}
\fi

\begin{document}

\title{The biHecke monoid of a finite Coxeter group}

\author{Florent Hivert\addressmark{1}, Anne Schilling\addressmark{2},
and Nicolas M.~Thiéry\addressmark{2,3}}

\address{\addressmark{1} LITIS (EA 4108), Université de Rouen,
  Avenue de l'Université  BP12
  76801 Saint-Etienne du Rouvray, France and
  Institut Gaspard Monge (UMR 8049)\\
\addressmark{2}Department of Mathematics, University of California, One
 Shields Avenue, Davis, CA 95616, U.S.A.\\
\addressmark{3}Univ Paris-Sud, Laboratoire de Mathématiques d'Orsay,
  Orsay, F-91405; CNRS, Orsay, F-91405, France}


\maketitle

\begin{abstract}
  \paragraph{Abstract.}
  The usual combinatorial model for the $0$-Hecke algebra
  $H_0(\sg[n])$ of the symmetric group is to consider the algebra (or
  monoid) generated by the bubble sort operators%
  . This construction generalizes  to any
  finite Coxeter group $W$. The authors previously introduced the
  Hecke group algebra, constructed as the algebra generated
  simultaneously by the bubble sort and antisort operators, and
  described its representation theory.
  In this paper, we consider instead the \emph{monoid} generated by
  these operators. We prove that it has $|W|$ simple and projective
  modules. In order to construct a combinatorial model for the simple
  modules, we introduce for each $w\in W$ a combinatorial module $T_w$
  whose support is the interval $[1,w]_R$ in right weak order.  This
  module yields an algebra, whose representation theory generalizes
  that of the Hecke group algebra. This involves the introduction of a
  $w$-analogue of the combinatorics of descents of $W$ and a
  generalization to finite Coxeter groups of blocks of permutation
  matrices.

    \paragraph{R\'esum\'e.}
  Le modèle combinatoire usuel pour la $0$-algèbre de Hecke
  $H_0(\sg[n])$ du groupe symétrique est obtenu en considérant
  l'algèbre (ou le monoïde) engendrée par les opérateurs de tri-à-bulle
  élémentaires,
  et cette construction se généralise pour tout groupe de Coxeter
  fini. Les auteurs ont précédemment introduit l'algèbre de Hecke
  groupe, engendrée conjointement par les opérateurs de tri et
  d'antitri à bulle.
  Dans cet article, nous considérons le monoïde engendré par ces
  opérateurs. Nous montrons qu'il a $|W|$ modules simples et
  projectifs. Afin de construire un modèle combinatoire pour ses
  modules simples, nous introduisons pour tout $w\in W$ un module
  combinatoire $T_w$ dont le support est l'intervalle $[1,w]_R$ pour
  l'ordre faible droit. Ce module détermine une algèbre dont la
  théorie des représentation généralise celle de l'algèbre de Hecke
  groupe, \emph{via} l'introduction d'un $w$-analogue de la
  combinatoire des descentes de $W$ et d'une généralisation aux
  groupes de Coxeter finis des blocs dans les matrices de
  permutations.  
\end{abstract}

\section{Introduction}

The usual combinatorial model for the $0$-Hecke algebra $H_0(\sg[n])$
of the symmetric group is the algebra (or monoid) generated by the
(anti) bubble sort operators $\pi_1,\ldots,\pi_{n-1}$, where $\pi_i$
acts on words of length $n$ and sorts the letters in positions $i$ and
$i+1$ decreasingly. By symmetry, one can also construct the bubble
sort operators $\opi_1,\ldots,\opi_{n-1}$, where $\opi_i$ acts by
sorting increasingly, and this gives an isomorphic construction
$\overline H_0$ of the $0$-Hecke algebra. This construction
generalizes naturally to any finite Coxeter group $W$.  Furthermore,
when $W$ is a Weyl group, and hence can be affinized, there is an
additional operator $\pi_0$ projecting along the highest root.

In~\cite{Hivert_Thiery.HeckeGroup.2007} the first and last author
constructed the \emph{Hecke group algebra} $\biheckealgebra$ by gluing
together the $0$-Hecke algebra and the group algebra of $W$ along
their right regular representation. Alternatively, $\biheckealgebra$
can be constructed as the \emph{biHecke algebra} of $W$, by gluing
together the two realizations $H_0(W)$ and $\overline H_0(W)$ of the
$0$-Hecke algebra. $\biheckealgebra$ admits a more conceptual
description as the algebra of all operators on $\K.W$ preserving left
antisymmetries; the representation theory of $\biheckealgebra$
follows, governed by the combinatorics of descents.
In~\cite{Hivert_Schilling_Thiery.HeckeGroupAffine.2008}, the authors
further proved that, when $W$ is a Weyl Group, $\biheckealgebra$ is a
natural quotient of the affine Hecke algebra.

In this paper, following a suggestion of Alain Lascoux, we study the
\emph{biHecke monoid} $\biheckemonoid(W)$, obtained by gluing together
the two $0$-Hecke \emph{monoids}. This involves the combinatorics of
the usual poset structures on $W$ (left, right, left-right, Bruhat
order), as well as a new one, the cutting poset, which in type $A$ is
related to blocks in permutation matrices. The guiding principle is
the use of representation theory to derive a (so far elusive)
summation formula for the size of this monoid, using that the simple
and projective modules of $\biheckemonoid$ are indexed by the elements
of $W$.

In type $A$, the tower of algebras $(\K[\biheckemonoid(\sg)])_{n\in
  \N}$ possesses long sought-after properties. Indeed, it is
well-known that several combinatorial Hopf algebras arise as
Grothendieck rings of towers of algebras. The prototypical example is
the tower of algebras of the symmetric groups which gives rise to the
Hopf algebra $\sym$ of symmetric functions, on the Schur
basis. Another example, due to Krob and
Thibon~\cite{Krob_Thibon.NCSF4.1997}, is the tower of the $0$-Hecke
algebras of the symmetric groups which gives rise to the Hopf algebra
$\qsym$ of quasi-symmetric functions of~\cite{Gessel.QSym.1984}, on
the $F_I$ basis.  The product rule on the $F_I$'s is naturally lifted
through the descent map to a product on permutations, leading to the
Hopf algebra $\fqsym$ of free quasi-symmetric functions.  This calls
for the existence of a tower of algebras $(A_n)_{n\in \N}$, such that
each $A_n$ contains $H_0(\sg)$ and has its simple modules indexed by
the elements of $\sg$. The biHecke monoids $M(\sg)$, and their Borel
submonoids $M_1(\sg)$, satisfy these properties, and are therefore
expected to yield new representation theoretical interpretations of
the bases of $\fqsym$.

In the remainder of this introduction, we briefly review Coxeter groups
and their $0$-Hecke monoids, and introduce our main objects of study:
the biHecke monoids.

\subsection{Coxeter groups}

Let $(W, S)$ be a Coxeter group, that is a group $W$ with a
presentation
\begin{equation}
  W = \langle\, S\ \suchthat\  (ss')^{m(s,s')},\ \forall s,s'\in S\,\rangle\,,
\end{equation}
with $m(s,s') \in \{1,2,\dots,\infty\}$ and $m(s,s)=1$.
The elements $s\in S$ are called \emph{simple reflections}, and the
relations can be rewritten as $s^2=\one$, where $\one$ is the identity in $W$ and:
\begin{equation}
    \underbrace{ss'ss's\cdots}_{m(s,s')} =
    \underbrace{s'ss'ss'\cdots}_{m(s,s')}\,, \qquad \text{ for all $s,s'\in S$}\, .
\end{equation}

Most of the time, we just write $W$ for $(W,S)$. In general, we follow
the notation of~\cite{Bjorner_Brenti.2005}, and we refer to this
monograph and to~\cite{Humphreys.1990} for details on Coxeter groups
and their Hecke algebras.  Unless stated otherwise, we always assume
that $W$ is finite, and denote its generators by $S=(s_i)_{i\in I}$
where $I=\{1,2,\ldots, n\}$ is the \emph{index set} of $W$.

The prototypical example is the Coxeter group of type $A_{n-1}$ which is the $n$-th symmetric
group $(W,S) := (\sg[n], \{s_1,\dots,s_{n-1}\})$, where
$s_i$ denotes the \emph{elementary transposition} which exchanges $i$ and $i+1$. 
The relations are $s_i^2=\one$ for $1\leq i \leq n-1$ and the \emph{braid relations}:
\begin{equation}
  \begin{alignedat}{2}
     s_i s_j         & = s_j s_i\,,            &     & \text{ for } |i-j|\geq2\,, \\
    s_i s_{i+1} s_i & = s_{i+1} s_i s_{i+1}\,, &\quad& \text{ for } 1\leq i\leq n-2\, .
  \end{alignedat}
\end{equation}
When writing a permutation $\mu \in \sg[n]$ explicitly, we use the
\emph{one-line notation}, that is the sequence
$\mu_1\mu_2\cdots\mu_n$, where $\mu_i:=\mu(i)$.

A \emph{reduced word} $i_1 \ldots i_k$ for an element $\mu \in W$
corresponds to a decomposition $\mu=s_{i_1}\cdots s_{i_k}$ of $\mu$
into a product of generators in $S$ of minimal length $k=\ell(\mu)$. A
\emph{(right) descent} of $\mu$ is an element $i\in I$ such that
$\ell(\mu s_i) < \ell(\mu)$. If $\mu$ is a permutation, this
translates into $\mu_i > \mu_{i+1}$. \emph{Left descents} are defined
analogously. The sets of left and right descents of $\mu$ are denoted
by $\Rec(\mu)$ and $\Des(\mu)$, respectively.

A Coxeter group $W$ comes equipped with four natural graded poset
structures. Namely $\mu<\nu$, in \emph{Bruhat order} (resp. \emph{left
  (weak)}, \emph{right (weak)}, \emph{left-right (weak) order}) if
some reduced word for $\mu$ is a subword (resp. right factor, left
factor, factor) of some reduced word for $\nu$.  In type $A$, the left
(resp. right) order give the usual left (resp. right) permutahedron.

For $J\subseteq I$, we denote by $W_J = \langle s_j \mid j \in J
\rangle$ the subgroup of $W$ generated by $s_j$ with $j\in
J$. Furthermore, the longest element in $W_J$ (resp. $W$) is denoted
by $s_J$ (resp. $w_0$).

\subsection{The $0$-Hecke monoid}
The \emph{$0$-Hecke monoid} $H_0(W) = \langle \pi_i \mid i \in I
\rangle$ of a Coxeter group $W$ is generated by the \emph{simple
  projections} $\pi_i$ with relations%
$\pi_i^2=\pi_i$ for $i\in I$ and
\begin{equation}
     \underbrace{\pi_i\pi_j\pi_i\pi_j\cdots}_{m(s_i,s_{i'})} =
    \underbrace{\pi_j\pi_i\pi_j\pi_i\cdots}_{m(s_i,s_{i'})} \qquad \text{ for all $i,j\in I$}\ .
\end{equation}
Thanks to these relations, the elements of $H_0(W)$ are canonically
indexed by the elements of $W$ by setting $\pi_w :=
\pi_{i_1}\cdots\pi_{i_k}$ for any reduced word $i_1 \dots i_k$ of $w$.
We further denote by $\pi_J$ the longest element of the
\emph{parabolic submonoid} $H_0(W_J) := \langle \pi_i \mid i \in J
\rangle$.

The right regular representation of $H_0(W)$ induces a concrete
realization of $H_0(W)$ as a monoid of operators acting on $W$, with
generators $\pi_1,\dots,\pi_n$ defined by:
\begin{equation}
  \label{equation.antisorting_action}
  w. \pi_i := \begin{cases}
    w & \text{if $i \in \Des(w)$,}\\
    ws_i & \text{otherwise.}
  \end{cases}
\end{equation}
In type $A$, $\pi_i$ sorts the letters at positions $i$ and $i+1$
decreasingly, and for any permutation $p$, $p.\pi_{w_0}=n\cdots21$.
This justifies naming $\pi_i$ an \emph{elementary bubble antisorting
  operator}.

Another concrete realization of $H_0(W)$ can be obtained by
considering instead the \emph{elementary bubble sorting operators}
$\opi_1,\dots,\opi_n$, whose action on $W$ are defined by:
\begin{equation}
  \label{equation.sorting_action}
  w. \opi_i := \begin{cases}
    ws_i  & \text{if $i\in \Des(w)$,}\\
    w & \text{otherwise.}
  \end{cases}
\end{equation}
In type $A$, and for any permutation $p$, one has $p.\opi_{w_0}=12\cdots n$.

\subsection{The biHecke monoid}

\begin{definition}
  Let $W$%
  be a finite Coxeter group. The \emph{biHecke monoid} is the
  submonoid of functions from $W$ to $W$ generated simultaneously by
  the elementary bubble sorting and antisorting operators
  of~\eqref{equation.antisorting_action}
  and~\eqref{equation.sorting_action}:
  \begin{equation*}
    \biheckemonoid := \biheckemonoid(W) :=
    \langle \pi_1, \pi_2, \ldots, \pi_n, \opi_1, \opi_2, \ldots, \opi_n \rangle\,.
  \end{equation*}
\end{definition}

\subsection{Outline}
The remainder of this paper is organized as follows.
In Section~\ref{section.blocks_cutting_poset}, we generalize the
notion of blocks of permutation matrices to any Coxeter group, and use
it to define a new poset structure on $W$, which we call the
\emph{cutting poset}; we prove that it is (almost) a lattice, and
derive that its M\"obius function is essentially that of the
hypercube.

In Section~\ref{section.M.combinatorics}, we study the combinatorial
properties of $\biheckemonoid(W)$. In particular, we prove that it
preserves left and Bruhat order, derive consequences on the fibers and
image set of its elements, prove that it is acyclic, and study its
conjugacy classes of idempotents.

In Section~\ref{section.M1}, our strategy is to consider a ``Borel''
triangular submonoid of $\biheckemonoid(W)$ whose representation
theory is simpler, but with the same number of simple modules, in the
hope to later induce back information about the representation theory
of $\biheckemonoid(W)$. Namely, we study the submonoid
$\biheckemonoid_1(W)$ of the elements fixing $\one$ in
$\biheckemonoid(W)$.  This monoid not only preserves Bruhat order, but
furthermore is \contracting.  It follows that it is $J$-trivial which
is the desired triangularity property. It is for example easily
derived that $\biheckemonoid_1(W)$ has $|W|$ simple modules, all of
dimension $1$. In fact most of our results about $\biheckemonoid_1$
generalize to any $J$-trivial monoid, which is the topic of a separate
paper on the representation theory of $J$-trivial
monoids~\cite{Denton_Hivert_Schilling_Thiery.JtrivialMonoids}.

In Section~\ref{section.translation_modules}, we introduce, for each
$w \in W$, the \emph{translation modules} $T_w$, whose support is the
interval $[1,w]_R$ in right order. Our original motivation, backed by
computer evidence, is that $T_w$ is closely related to the
indecomposable projective module $P_w$ of $\biheckemonoid(W)$. In
particular it is indecomposable, and we can use $T_w$ to construct a
combinatorial model for the simple module $S_w$ of $\biheckemonoid(W)$
which appears as the top of $T_w$. We derive a formula for the
dimension of $S_w$, using an inclusion-exclusion on the sizes of
intervals in $(W, \le_R)$, along the cutting poset. On the way, we
study the algebra $\wbiheckealgebra{w}$ induced by the action of
$\biheckemonoid(W)$ on $T_w$. It turns out to be a natural
$w$-analogue of the Hecke group algebra, acting not anymore on the
full Coxeter group, but on the interval $[1,w]_R$ in right order.  All
the properties of the Hecke group algebra pass through this
generalization, with the combinatorics of descents being replaced by
that of blocks and of the cutting poset. In particular,
$\wbiheckealgebra{w}$ is Morita equivalent to the incidence algebra of
a lattice.  

In Section~\ref{section.M.representation_theory}, we derive (parts of)
the representation theory of $\biheckemonoid(W)$ from
Sections~~\ref{section.M.combinatorics}, \ref{section.M1},
and~\ref{section.translation_modules}.

A long version of this paper with all proofs included will appear separately.

\subsection*{Acknowledgments}
We would like to thank Jason Bandlow, Mahir Can, Brant Jones, Alain
Lascoux, Jean-Christophe Novelli, Jean-Éric Pin, Vic Reiner, Franco
Saliola, Benjamin Steinberg, and Jean-Yves Thibon for enlightening
discussions.
This research was driven by computer exploration, using the
open-source mathematical software \texttt{Sage}~\cite{sage} and its
algebraic combinatorics features developed by the
\texttt{Sage-Combinat} community~\cite{Sage-Combinat}.  AS and NT were in part supported by NSF grants DMS--0652641 and DMS--0652652.
AS was in part supported by the Kastler foundation.
FH was partly supported by ANR grant 06-BLAN-0380.

\section{Blocks of Coxeter group elements and the cutting poset}
\label{section.blocks_cutting_poset}

In this section, we develop the combinatorics underlying the
representation theory of the translation modules introduced in
Section~\ref{section.translation_modules}. The key question is: for
which subsets $J \subseteq I$ does the canonical bijection between a Coxeter group
$W$ and the Cartesian product $W_J \times \lcoset{W}{J}$ of a
parabolic subgroup $W_J$ by its set of cosets representatives
$\lcoset{W}{J}$ in $W$ restrict properly to an interval $[1,w]_R$ in
right order (see Figure~\ref{figure.4312})? In type $A$, the answer is
given by the so-called \emph{blocks} in the permutation matrix of $w$,
and we generalize this notion to any Coxeter group.
After reviewing parabolic subgroups and cosets representatives in
Section~\ref{section.parabolic}, we define blocks of Coxeter group
elements in Section~\ref{section.blocks} and show in
Section~\ref{section.tiling.type_A} how this notion specializes to
type $A$.  Finally, in Section~\ref{section.cutting_poset}, we
introduce and study the cutting poset.%

\subsection{Parabolic subgroups and cosets representatives}
\label{section.parabolic}

For a subset $J \subseteq I$, the \emph{parabolic subgroup} $W_J$ of $W$ is
the Coxeter subgroup of $W$ generated by $s_j,j\in J$.  A complete
system of minimal length representatives of the right cosets $W_J w$
(resp. of the left cosets $w W_J$) are given respectively by:
\begin{equation*}
    \lcoset{W}{J} := \{ x \in W \suchthat \Rec(x) \cap J = \emptyset\} \qquad \text{and} \qquad
    \rcoset{W}{J} := \{ x \in W \suchthat \Des(x) \cap J = \emptyset\}\,.
\end{equation*}

Every $w\in W$ has a unique decomposition $w = w_J \lcoset{w}{J}$ with
$w_J\in W_J$ and $\lcoset{w}{J} \in \lcoset{W}{J}$. Similarly, there
is a unique decomposition $w = \rcoset{w}{K} {}_Kw$ with ${}_Kw \in
{}_KW=W_K$ and $\rcoset{w}{K} \in \rcoset{W}{K}$.
A subset $J$ is \emph{left reduced} w.r.t. $w$ if $J'\subset J$
implies $\lcoset{w}{J'}>_L \lcoset{w}{J}$. Right reduced $K$'s are
defined analogously.

\subsection{Blocks of Coxeter group elements}
\label{section.blocks}

We now come to the definition of blocks of Coxeter group elements, and
associated cutting points.

\begin{definition}[Blocks and cutting points]
  \label{definition.block}
  \label{definition.cutting_point}
  Let $w\in W$.  We say $K\subseteq I$ is a \emph{right block}
  (resp. $J \subseteq I$ is a \emph{left block}) of $w$, if there
  exists $J\subseteq I$ (resp. $K\subseteq I$) such that
    $w W_K = W_J w\,.$

  In that case, $v:=w^K$ is called a \emph{cutting point} of $w$,
  which we denote by $v \cuteq w$. Furthermore, $K$ is \emph{proper}
  if $K\neq \emptyset$ and $K\neq I$; it is \emph{nontrivial} if
  $\rcoset{w}{K}\neq w$ (or equivalently ${}_K w\neq 1$); analogous
  definitions are made for left blocks.

  We denote by $\Ks{w}$ (resp. $\Js{w}$) the set of all right
  (resp. left) blocks for $w$, and by $\RKs{w}$ (resp. $\RJs{w}$) the
    set of all reduced right (resp. left) blocks for $w$.
  
\end{definition}

\begin{proposition} \label{proposition.tiling}
    Assuming that $J$ is reduced, $J$ is a left block of $w$ if and only if the bijection
  \begin{equation*}
    \begin{cases}
      W_J \times \lcoset{W}{J} &\to W\\
      (u,v) &\mapsto uv
    \end{cases}
  \end{equation*}
  restricts to a bijection $[1,w_J]_R \times [1,\lcoset{w}{J}]_R \to
  [1,w]_R$ (see Figure~\ref{figure.4312}).
\end{proposition}

Due to Proposition~\ref{proposition.tiling}, we also say that $[1,v]_R$
\emph{tiles} $[1,w]_R$ if $v=\lcoset{w}{J}$ for some left block $J$%
.

\begin{example}
  For $w=w_0$, any $K\subseteq I$ is a reduced right block; of course
  $\;\rcoset{w_0}{K} \le_L w_0$, and ${}_K w_0$ is the maximal element
  of the parabolic subgroup $W_K = {}_K W$. The cutting point
  $w^K\cuteq w$ is the maximal element of the right descent class for
  the complement of $K$.

\end{example}

\begin{proposition}
  \label{proposition.block.union}
  The set $\Js{w}$ (resp. $\Ks{w}$) of left (resp. right) blocks is
  stable under union and intersection. So, they form a sublattice of
  the Boolean lattice.

  The sets $\RJs{w}$ and $\RKs{w}$ are (dual) Moore families and
  therefore lattices.
\end{proposition}

\begin{definition}[$w$-codescent sets]
  \label{definition.max_block}
  For $u\in [1,w]_R$ define $\Kblock{w}(u)$ to be the maximal reduced
  right block $K$ of $w$ such that $u$ is below the corresponding
  cutting point, that is $u\le_R w^K$.  Let $\Jblock{w}(u)$ be the
  left block corresponding to $\Kblock{w}(u)$.
\end{definition}

\begin{example}
  When $w=w_0$, any $J\subseteq I$ is a reduced left
  block. Furthermore, for $u\in W$, $\Jblock{w_0}(u)$ is the
  complement of its left-descent set: %
      $\Jblock{w_0}(u) = I \setminus \Rec(u)$. Idem on the right.

\end{example}

\subsection{Blocks of permutations}
\label{section.tiling.type_A}
We now specialize to type $A_{n-1}$. For a permutation $w \in \sg[n]$,
the blocks of Definition~\ref{definition.block} correspond to the
usual notion of blocks of the matrix of $w$ (or unions thereof), and
the cutting points $w^K$ for right blocks $K$ correspond to putting
the identity in the matrix-blocks of $w$.
A matrix-block of a permutation $w$ is an interval $[k',k'+1,\dots,k]$
which is mapped to another interval. Pictorially, this corresponds to
a square submatrix of the matrix of $w$ which is again a permutation
matrix (that of the \emph{associated permutation}).  For example, the
interval $[2,3,4,5]$ is mapped to the interval $[4,5,6,7]$ by the
permutation $w=36475812\in \sg[8]$, and is therefore a matrix-block of
$w$ with associated permutation $3142$. Similarly, $[7,8]$ is a
matrix-block with associated permutation $12$:
\begin{center}
\scalebox{.6}[-.6]{
$
  \begin{array}{|c|cccc|c|cc|}
    \hline
     & & & & &\bullet &&\\
    \hline
     & & &\bullet& & & &\\
     &\bullet& & & & & &\\
     & & & &\bullet& & &\\
     & &\bullet& & & & &\\
    \hline
    \bullet& & & & & & &\\
    \hline
     & & & & & & &\bullet\\
     & & & & & &\bullet &\\
    \hline
  \end{array}
$}
\end{center}

For any permutation $w$, the singletons $[i]$ and the full set
$[1,2,\dots,n]$ are always matrix-blocks; the other matrix-blocks of
$w$ are called \emph{proper}.  A permutation with no proper
matrix-block, such as $58317462$, is called \emph{simple}.
See~\cite{Albert_Atkinson.2005} for a review of simple permutations.

\begin{proposition}
  \label{proposition.blocks.matrix_blocks}
  Let $w\in \sg[n]$. The right blocks of $w$ are in bijection with
  disjoint unions of (non singleton) matrix-blocks for $w$; each
  matrix-block with column set $[i,i+1,\dots,k]$ contributes
  $\{i,i+1,\dots,k-1\}$ to the right block; each matrix-block with row
  set $[i,i+1,\ldots,k]$ contributes $\{i,i+1,\dots,k-1\}$ to the left
  block. %
  In addition, trivial right blocks correspond to unions of identity
  matrix-blocks.  Also, reduced right blocks correspond to unions of
  connected matrix-blocks.
\end{proposition}

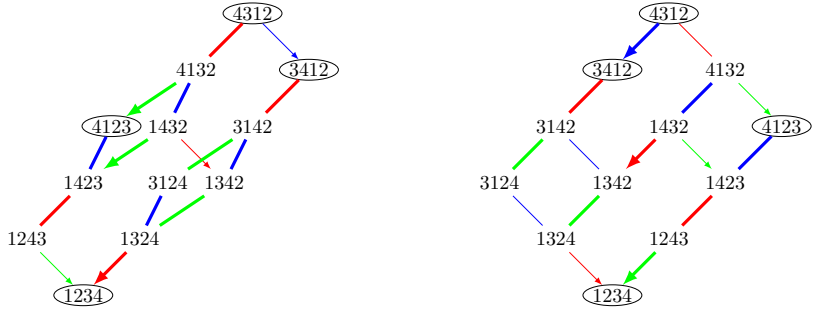
\begin{figure}
  
    \centerline{
    \scalebox{.75}{\tikzstyle{cutting}=[shape=ellipse,draw=black,inner sep=.5pt]
\tikzstyle{interval}=[ultra thick]
\begin{tikzpicture}[baseline=(current bounding box.east), >= latex]
  \node[cutting] (1234) at ( 0, 0) {1234};

  \node (1324) at ( 1, 1) {1324};

  \node (3124) at ( 1.5, 2) {3124};
  \node (1342) at ( 2.5, 2) {1342};

  \node (3142) at ( 3, 3) {3142};
  \node[cutting] (3412) at ( 4, 4) {3412};

  \node (1243) at (-1, 1) {1243};
  \node (1423) at ( 0, 2) {1423};

  \node[cutting] (4123) at (0.5, 3) {4123};
  \node (1432) at ( 1.5, 3) {1432};

  \node (4132) at (2, 4) {4132};

  \node[cutting] (4312) at ( 3, 5) {4312};

  \draw[<-,red,interval]   (1234)  -- (1324);
  \draw[<-,green] (1234)  -- (1243);

  \draw[blue,interval]  (1324)  -- (3124);
  \draw[green,interval] (1324)  -- (1342);

  \draw[red,interval]   (1243)  -- (1423);

  \draw[<-,red]   (1342)  -- (1432);
  \draw[green,interval] (3124)  -- (3142);
  \draw[blue,interval]  (1342)  -- (3142);

  \draw[blue,interval]  (1423)  -- (4123);
  \draw[<-,green,interval] (1423)  -- (1432);

  \draw[red,interval]   (3142)  -- (3412);

  \draw[blue,interval]  (1432)  -- (4132);
  \draw[<-,green,interval] (4123)  -- (4132);

  \draw[<-,blue]  (3412)  -- (4312);
  \draw[red,interval]   (4132)  -- (4312);
\end{tikzpicture}}
    \qquad \qquad
    \scalebox{.75}{\tikzstyle{cutting}=[shape=ellipse,draw=black,inner sep=.5pt]
\tikzstyle{interval}=[ultra thick]
\begin{tikzpicture}[baseline=(current bounding box.east), >=latex]
  \node[cutting] (1234) at ( 0, 0) {1234};

  \node (1324) at (-1, 1) {1324};
  \node (1243) at ( 1, 1) {1243};

  \node (3124) at (-2, 2) {3124};
  \node (1342) at ( 0, 2) {1342};
  \node (1423) at ( 2, 2) {1423};

  \node (3142) at (-1, 3) {3142};
  \node (1432) at ( 1, 3) {1432};
  \node[cutting] (4123) at ( 3, 3) {4123};

  \node[cutting] (3412) at ( 0, 4) {3412};
  \node (4132) at ( 2, 4) {4132};

  \node[cutting] (4312) at ( 1, 5) {4312};

  \draw[<-,red]   (1234)  -- (1324);
  \draw[<-,green,interval] (1234)  -- (1243);

  \draw[blue]  (1324)  -- (3124);
  \draw[green,interval] (1324)  -- (1342);

  \draw[red,interval]   (1243)  -- (1423);

  \draw[green,interval] (3124)  -- (3142);
  \draw[blue]  (1342)  -- (3142);
  \draw[<-,red,interval]   (1342)  -- (1432);

  \draw[blue,interval]  (1423)  -- (4123);
  \draw[<-,green] (1423)  -- (1432);

  \draw[red,interval]   (3142)  -- (3412);

  \draw[blue,interval]  (1432)  -- (4132);
  \draw[<-,green] (4123)  -- (4132);

  \draw[<-,blue,interval]  (3412)  -- (4312);
  \draw[red]   (4132)  -- (4312);
\end{tikzpicture}}
  }
  
  \caption{Two pictures of the interval $[1234,4312]_R$ in right order
    illustrating its proper tilings, for $J:=\{3\}$ and $J:=\{1,2\}$,
    respectively. The thick edges highlight the tiling. The circled
    permutations are the cutting points, which are at the top of the
    tiling intervals.  Blue, red, green lines correspond to $s_1$,
    $s_2$, $s_3$ respectively. See Section~\ref{subsection.bihecke}
    for the definition of the orientation of the edges (this is
    $G^{(4312)}$); edges with no arrow tips point in both directions.}
  \label{figure.4312}
\end{figure}

\begin{example}
  As in Figure~\ref{figure.4312}, consider the permutation $4312$, whose permutation matrix is:
  \begin{center}
  \scalebox{.6}[-.6]{$
    \begin{array}{|c|c|cc|}
      \hline
      \bullet&&&\\\hline
      &\bullet&&\\    \hline
      & & &\bullet\\
      & &\bullet &\\
      \hline
    \end{array}
  $}
  \end{center}
  The reduced (right)-blocks are $K=\{\}$, $\{1\}$, $\{2,3\}$, and
  $\{1,2,3\}$. The cutting points are $4312$, $3412$, $4123$, and
  $1234$, respectively. The corresponding left blocks are $J=\{\}$,
  $\{3\}$, $\{1,2\}$ and $\{1,2,3\}$, respectively.
  The non-reduced blocks are $\{3\}$ and $\{1,3\}$, as they are
  respectively equivalent to the blocks $\{\}$ and $\{1\}$. Finally,
  the trivial blocks are $\{\}$ and $\{3\}$.
\end{example}

\subsection{The cutting poset}
\label{section.cutting_poset}

\begin{theorem}
  \label{theorem.cutting}
  $(W,\cuteq)$ is a poset with $\one$ as minimal element; it is further a
  subposet of both left and right order.
  Every interval of $(W,\cuteq)$ is a sublattice of both left and
  right order.

  The $\cuteq$-lower covers of an element $w$ correspond to the
  nontrivial blocks of $w$ which are minimal for inclusion. The
  meet-semilattice $L_w$ they generate is free, and is in
  correspondence with the lattice of unions of these minimal
  nontrivial blocks, or alternatively of the intersections of the
  intervals $[1,u]_R$ for $u$ $\cuteq$-lower covers of $w$.

  The M\"obius function is given by $\mu(u,w)=\pm 1$ if $u$ is in $L_w$
  (with alternating sign according to the usual rule for the Boolean lattice), and $0$
  otherwise.
\end{theorem}
This M\"obius function is used
in %
Section~\ref{subsection.bihecke} to compute the size of the simple
modules of $\biheckemonoid$.
\begin{conjecture}
  $(W,\cuteq)$ is a meet-semilattice whose intervals are all
  distributive lattices.
\end{conjecture}

\section{The combinatorics of $\biheckemonoid(W)$}
\label{section.M.combinatorics}

In this section we study the combinatorics of the biHecke monoid
$\biheckemonoid(W)$ of a finite Coxeter group $W$. In particular, we
prove in Section~\ref{ss.left} that its elements
preserve left order and Bruhat order, and derive in
Section~\ref{ss.fiber_image} properties of their image sets and
fibers.
Finally, in Section~\ref{ss.idempotents},
we prove the key combinatorial ingredients
for the enumeration of the simple modules of $\biheckemonoid(W)$ in
Section~\ref{section.M.representation_theory}: $\biheckemonoid(W)$ is
acyclic and admits $|W|$ conjugacy classes of idempotents.

\subsection{Preservation of left and Bruhat order}
\label{ss.left}

\begin{lemma}

  Take $f\in \biheckemonoid(W)$, $w\in W$, and $j\in I$. Then, $(s_j w).f$ is either $w.f$ or
  $s_j(w.f)$.
\end{lemma}

\begin{proposition}
  \label{proposition.weak_order2}
    The elements $f$ of $\biheckemonoid$ preserve left order:
  $u\leq_L v\Rightarrow u.f\leq_L v.f$.
  
\end{proposition}

\begin{proposition}
  \label{proposition.bruhat}
  The elements $f$ of $\biheckemonoid$ preserve Bruhat order:
  $u\leq_B v\Rightarrow u.f\leq_B v.f$.
\end{proposition}

\begin{proposition}
  \label{corollary.BruhatDecrease}
  Any $f\in\biheckemonoid$ such that $\one.f=\one$ is \contracting
  for Bruhat order: $w.f\le_B w$.
\end{proposition}

\subsection{Fibers and image sets}
\label{ss.fiber_image}

\begin{proposition}
  \label{proposition.image_set.idempotent}
  The image set of an idempotent in $\biheckemonoid(W)$ is an interval in left order.
\end{proposition}

\begin{proposition}
  \label{proposition.contract_fibers}
  Take $f\in \biheckemonoid(W)$, and consider the Hasse diagram of
  left order contracted with respect to the fibers of $f$. Then, this
  graph is isomorphic to left order restricted on the image set.
\end{proposition}

\begin{proposition}
  \label{proposition.fibers_image_set}
  Any element $f\in \biheckemonoid(W)$ is characterized by its set of
  fibers and its image set.
\end{proposition}

A monoid $M$ is called \emph{acyclic} if for any $f\in M$, there
exists $k>0$ such that $f^{k+1}=f^k$. Note that, in that case,
$f^\infty := f^k=f^{k+1}=\dots$ is an idempotent.

\begin{proposition}
  \label{proposition.pseudo_idempotent}
  The biHecke monoid $\biheckemonoid(W)$ is \emph{acyclic}.
\end{proposition}

\subsection{Conjugacy classes of idempotents}
\label{ss.idempotents}

\begin{proposition}\mbox{}
  \label{proposition.uniqueness_fix1}
   \label{e.uniqueness}
    For $w\in W$,
     $e_w := \pi_{w^{-1} w_0} \opi_{w_0 w}$
    is the unique idempotent with image set $[1,w]_L$.
    For $u\in W$, it satisfies $e_w(u) = \max_{\le_B} \bigl( [1,u]_B \cap [1,w]_L \bigr)$.

\end{proposition}

\begin{corollary}
  \label{corollary.e}
  For $u,w \in W$, the intersection $[1,u]_B \cap [1,w]_L$ is a lower
  $\le_L$ ideal with a unique maximal element $v$ in Bruhat order. The
  maximum is given by $v=e_w(u)$.
\end{corollary}

We
are now in the position to describe the conjugacy relations between the 
idempotents of $M$.

\begin{lemma}
  \label{lemma.conjugate}
  Let $e$ and $f$ be idempotents with respective image sets $[a,b]_L$
  and $[c,d]_L$. Then, $f\in MeM$ if and only if $dc^{-1} \le_{LR}
  ba^{-1}$.
  In particular, $e$ and $f$ are conjugates if and only if the
  intervals $[a,b]_L$ and $[c,d]_L$ are of the same type: $dc^{-1} =
  ba^{-1}$.
\end{lemma}

\begin{corollary}
  \label{corollary.transversal}
  The idempotents $(e_w)_{w\in W}$ form a complete set of
  representatives of the conjugacy classes of idempotents in
  $\biheckemonoid$.
\end{corollary}

\section{The Borel submonoid $\biheckemonoid_1(W)$ and its representation theory}
\label{section.M1}

In the previous section, we outlined the importance of the idempotents
$(e_w)_{w\in W}$. An crucial feature is that they live in a ``Borel''
submonoid %
$\biheckemonoid_1 := \{ f\in \biheckemonoid \mid \one.f = \one\}$.
In fact:
\begin{theorem}
  \label{theorem.M1.generating_set}
  $\biheckemonoid_1$ has a unique minimal generating set which
  consists of the ($2^n-n$ in type $A$) idempotents $e_w$ where $w_0
  w^{-1}$ is Grassmanian.
\end{theorem}

Furthermore, the elements of $\biheckemonoid_1$ are both order-preserving and \contracting for Bruhat order. %
\begin{corollary}
  \label{remark.M1}
  For $f,g \in \biheckemonoid_1$, define the relation $f\leq g$ if
  $w.f \leq w.g$ for all $w\in W$. Then, $\le$ defines a partial order
  on $\biheckemonoid_1$ such that $fg \leq f$ and $fg \leq g$ for any
  $f,g \in \biheckemonoid_1$.
\end{corollary}
A monoid with such an order is called \emph{$J$-trivial}, and the
description of its representation theory is the topic of a separate
paper~\cite{Denton_Hivert_Schilling_Thiery.JtrivialMonoids}. We
summarize here the main results for $\biheckemonoid_1$.

For each $w\in W$ define $S_w$ to be the one-dimensional vector space
with basis $\{\epsilon_w\}$ together with the right operation of any
$f\in \biheckemonoid_1$ given by %
  $\epsilon_w.f :=\epsilon_w$ if $w.f=w$ and $\epsilon_w.f:=0$ otherwise.
The basic features of the representation theory of $\biheckemonoid_1$
can be stated as follows:
\begin{theorem}
  The radical of $\K[\biheckemonoid_1]$ is the ideal with basis
  $(f^\infty-f)_f$, for $f$ non-idempotent.  The quotient of
  $\K[\biheckemonoid_1]$ by its radical is commutative. Therefore, all
  simple $\biheckemonoid_1$-module are one dimensional. In fact, the
  family $\{S_w\}_{w\in W}$ forms a complete system of representatives
  of the simple $\biheckemonoid_1$-modules.
\end{theorem}
To describe the indecomposable projective modules, we note that the
restriction of the conjugacy relation ($J$-order) to idempotents has a
very simple description: %
\begin{proposition}
  For $u,v\in W$, the following are equivalent:
  \begin{align*}
  \bullet\ e_u e_v &= e_u\,;&&
  \bullet\ \text{$v \leq_L u$ for left order;} \\
  \bullet\ e_v e_u &= e_u\,; &&
  \bullet\ \text{there exists $x,y\in \biheckemonoid_1$ such that $e_u = x e_v y$\,.}
  \end{align*}
  Moreover $(e_ue_v)^\infty=e_{u\join_L v}$, where $u\join_L v$ is the
    join  of $u$ and $v$ in
    left order.
\end{proposition}

\begin{definition}
  For any element $x\in M$, define 
  \begin{equation}
    \label{equation.lfix_rfix}
    \lfix(x) := \min_{\le_L} \{u\in W \suchthat e_ux=x\} \qquad \text{and} \qquad
    \rfix(x) := \min_{\le_L} \{u\in W \suchthat xe_u=x\}\,.
  \end{equation}
  
\end{definition}

Then, the projective modules and Cartan invariants can be described as follows:
\begin{theorem}
  There is an explicit basis $(b_f)_{f\in\biheckemonoid_1}$ of
  $\K[\biheckemonoid_1]$ such that, for all $w\in W$,
  \begin{itemize}
  \item the family $\{b_x \suchthat \lfix(x) = w\}$ is a basis for the right
    projective module associated to $S_w$;
  \item the family $\{b_x \suchthat \rfix(x) = w\}$ is a basis for the left
    projective module associated to $S_w$.
  \end{itemize}
  Moreover, the Cartan invariant of $\K[\biheckemonoid_1]$ defined by $c_{u,v} :=
  \dim(e_u \K[\biheckemonoid_1]e_v)$ for $u,v \in W$ is given by $c_{u,v} =
  |C_{u,v}|$, where
    $C_{u,v} := \{f\in \biheckemonoid_1 \suchthat u=\lfix(f)\text{ and } v=\rfix(f) \}$.

\end{theorem}

\section{Translation modules and $w$-biHecke algebras}
\label{section.translation_modules}

The main purpose of this section is to pave the ground for the
construction of the simple modules $S_w$ of the biHecke monoid
$\biheckemonoid=\biheckemonoid(W)$ in
Section~\ref{subsection.M.simple_modules}.  To this end, in
Section~\ref{subsection.translation}, we endow the interval $[1,w]_R$
with a natural structure of combinatorial $\biheckemonoid$-module
$T_w$, called \emph{translation module}. This module is closely
related to the projective module $P_w$ of $\biheckemonoid$
(Corollary~\ref{corollary.translation_module}), which explains its
important role.
By taking the quotient of $\K[\biheckemonoid]$ through its
representation on $T_w$, we obtain a $w$-analogue
$\wbiheckealgebra{w}$ of the biHecke algebra $\biheckealgebra$. This
algebra turns out to be interesting in its own right, and we proceed
by generalizing most of the results
of~\cite{Hivert_Thiery.HeckeGroup.2007} on the representation theory
of $\biheckealgebra$.

As a first step, we introduce in
Section~\ref{subsection.translation.P} a collection of submodules
$P_J^{(w)}$ of $T_w$, which are analogues of the projective modules of
$\biheckealgebra$. Unlike for $\biheckealgebra$, not any subset $J$ of
$I$ yields such a submodule, and this is where the combinatorics of
the blocks of $w$ enters the game. In a second step, we derive in
Section~\ref{subsection.translation.triangular} a lower bound on the
dimension of $\wbiheckealgebra{w}$; this requires a (fairly involved)
combinatorial construction of a family of functions on $[1,w]_R$ which
is triangular with respect to Bruhat order.
In Section~\ref{subsection.bihecke} we combine these results to derive
the dimension and representation theory of $\wbiheckealgebra{w}$:
projective and simple modules, Cartan matrix, quiver, etc.

\subsection{Translation modules and $w$-biHecke algebras}
\label{subsection.translation}
For $f\in M$,  define the \emph{type} of $f$ by
$\type(f) := (w_0.f)(\one.f)^{-1}$. %
By Proposition~\ref{proposition.weak_order2},
we know that for $f,g\in M$ either $\type(fg)=\type(f)$, or
$\ell(w_0.(fg))-\ell(\one.(fg)) < \ell(w_0.f) - \ell(\one.f)$ and
hence $\type(fg)\neq \type(f)$.  The second case occurs precisely when
$\fiber(f)$ is strictly finer than $\fiber(fg)$ or equivalently
$\rank(fg)<\rank(f)$, where the \emph{rank} is the cardinality of the image set.

\begin{definition}
  Fix $f\in M$. The right $M$-module
  \begin{equation*}
    \trans(f) := \K.fM/\K.\{h \in fM \mid \rank(h) < \rank(f)\}\\
  \end{equation*}
  is called the \emph{translation module} associated with $f$.
\end{definition}

\begin{proposition}
  Fix $f\in M$. Then:
  \begin{equation}
    \{h\in fM \suchthat \rank(h)=\rank(f)\} = \{ f_u \suchthat u\in [1,\type(f)^{-1}w_0]_R \}\,,
  \end{equation}
  where $f_u$ is the unique element of $M$ such that
  $\fiber(f_u)=\fiber(f)$ and $\one.f_u=u$.
\end{proposition}

\begin{proposition} \label{proposition.combinatorial_translation}
  Set $w=\type(f)^{-1}w_0$. Then, $(f_u)_{u\in [1,w]_R}$ forms a basis
  of $\trans(f)$ such that: %
  \def\unindent{\hspace{-.3cm}}
  \begin{equation}
    f_u.\pi_i =
    \begin{cases}
      f_u      &\unindent\text{if $i\in \Des(u)$}\\
      f_{u s_i} &\unindent\text{if $i\not\in \Des(u)$ and $us_i\in [1,w]_R$}\\
      0     & \unindent\text{otherwise;}
    \end{cases} \ 
    f_u.\opi_i =
    \begin{cases}
      f_{u s_i} & \unindent\text{if $i\in \Des(u)$ and $us_i\in [1,w]_R$}\\
      f_u      & \unindent\text{if $i\not\in \Des(u)$}\\
      0        & \unindent\text{otherwise;}
    \end{cases}
\end{equation}

\end{proposition}

This proposition gives a combinatorial model for translation
modules. It is clear that two functions with the same type yield
isomorphic translation modules. The converse also holds:
\begin{proposition}
  For any $f,f'\in M$, the translation modules $\trans(f)$ and $\trans(f')$ are isomorphic
  if and only if $\type(f) = \type(f')$.
\end{proposition}

By the previous proposition, we may choose a canonical representative
for translation modules. We choose $T_w := \trans(e_{w,w_0})$, and
identify its basis with $[1,w]_R$ via $u\mapsto f_u$.

\begin{definition}
  The $w$-biHecke algebra $\wbiheckealgebra{w}$ is the natural
  quotient of $\K[\biheckemonoid(W)]$ through its representation on
  $T_w$.  In other words, it is the subalgebra of $\End(T_w)$
  generated by the operators $\pi_i$ and $\opi_i$ of
  Proposition~\ref{proposition.combinatorial_translation}.
\end{definition}

\subsection{Left antisymmetric submodules}
\label{subsection.translation.P}

By analogy with the simple reflections in the Hecke group algebra, we
define for each $i\in I$ the operator $s_i := \pi_i + \opi_i -1$.
For $u\in [1,w]_R$, it satisfies
$u.s_i = us_i$ if $us_i\in [1,w]_R$ and $u.s_i = -u$ otherwise.
These operators are still involutions, but do not quite satisfy the
braid relations. %
One can further define operators $\ls_i$ acting similarly on the left.

\begin{definition}
  For $J\subseteq I$, set $P^{(w)}_J := \{ v\in T_w \suchthat \ls_i.v=-v, \quad
  \forall i \in J\}$.
\end{definition}
For $w=w_0$, these are the projective modules $P_J$ of the biHecke
algebra.

\begin{proposition}
  \label{proposition.cond_Jtiling}
  Take $w\in W$ and $J\subseteq I$ left reduced%
  . Then,
    $J$ is a reduced left block of $w$ if and only if $P^{(w)}_J$ is an $M$-submodule of $T_w$.
  
\end{proposition}

  It is clear from the definition that for $J_1, J_2 \subseteq I$,
  $P^{(w)}_{J_1\cup J_2} = P^{(w)}_{J_1} \cap P^{(w)}_{J_2}$.
  Since the set $\Js{w}$ of reduced left blocks of $w$ is stable under
  union, the set of $M$-modules $(P^{(w)}_J)_{J\in \Js{w}}$ is stable
  under intersection.
On the other hand, unless $J_1$ and $J_2$ are comparable,
$P^{(w)}_{J_1\cup J_2}$ is a strict subspace of $P^{(w)}_{J_1} +
P^{(w)}_{J_2}$. Hence, for $J\in \Js{w}$, we set $S^{(w)}_J:=
P^{(w)}_J / \sum_{J'\supsetneq J, J' \in \Js{w}} P^{(w)}_{J'}$.

\subsection{A (maximal) Bruhat-triangular family of $\wbiheckealgebra{w}$}
\label{subsection.translation.triangular}

Consider the submonoid $F$ in $\wbiheckealgebra{w}$ generated by the
operators $\pi_i$, $\opi_i$, and $s_i$, for $i\in I$. For $f\in F$ and
$u\in[1,w]_R$, we have $u.f = ± v$ for some $v\in[1,w]_R$. For our
purposes, the signs can be ignored and $f$ be considered as a function
from $[1,w]_R$ to $[1,w]_R$.

\begin{definition}
  For $u,v \in [1,w]_R$, a function $f\in F$ is called
  \emph{$(u,v)$-triangular} (for Bruhat order) if $v$ is the unique
  minimal element of $\im(f)$ and $u$ is the unique maximal element of
  $f^{-1}(v)$ (all minimal and maximal elements in this context are
  with respect to Bruhat order).
\end{definition}

\begin{proposition}
  \label{proposition.triangular}
  Take $u,v \in [1,w]_R$ such $\Kblock{w}(u) \subseteq
  \Kblock{w}(v)$. Then, there exists a $(u,v)$-triangular function
  $f_{u,v}$ in $F$.
\end{proposition}
For example, for $w=4312$ in $\sg[4]$, the condition on $u$ and $v$ is
equivalent to the existence of a path from $u$ to $v$ in the digraph
$G^{(4312)}$ (see Figure~\ref{figure.4312} and Section~\ref{subsection.bihecke}).

The construction of $f_{u,v}$ is explicit, and the triangularity
derives from $f_{u,v}$ being either in $M$, or close enough to be
bounded below by an element of $M$.  It follows from the upcoming
Theorem~\ref{theorem.wBihecke.representations} that the condition on
$u$ and $v$ is not only sufficient but also necessary.

\subsection{Representation theory of $w$-biHecke algebras}
\label{subsection.bihecke}

Consider the digraph $G^{(w)}$ on $[1,w]_R$ with an edge $u\mapsto v$
if $u=vs_i$ for some $i$, and $\Jblock{w}(u) \subseteq
\Jblock{w}(v)$. Up to orientation, this is the Hasse diagram of right
order (see for example Figure~\ref{figure.4312}).
The following theorem is a generalization of~\cite[Section 3.3]{Hivert_Thiery.HeckeGroup.2007}.
\begin{theorem}
  \label{theorem.wBihecke.representations}
  $\wbiheckealgebra{w}$ is the maximal algebra stabilizing all the
  modules $P_J^{(w)}$, for $J\in \RJs{w}$%
  .

  The elements $f_{u,v}$ of Proposition~\ref{proposition.triangular}
  form a basis $\wbiheckealgebra{w}$; in particular,
  \begin{equation}
    \dim \wbiheckealgebra{w} = |\{ (u,v) \suchthat \Jblock{w}(u)
    \subseteq \Jblock{w}(v) \}| \; .
  \end{equation}

  $\wbiheckealgebra{w}$ is the digraph algebra of the graph $G^{(w)}$.

  The family $(P_J)_{J\in \Js{w}}$ forms a complete system of
  representatives of the indecomposable projective modules of
  $\wbiheckealgebra{w}$.

  The family $(S_J)_{J\in \Js{w}}$ forms a complete system of
  representatives of the simple modules of $\wbiheckealgebra{w}$. The
  dimension of $S_J$ is the size of the corresponding $w$-descent class.

  $\wbiheckealgebra{w}$ is Morita equivalent to the poset algebra of
  the lattice $[1,w]_\cuteq$.
\end{theorem}

\section{The representation theory of $\biheckemonoid(W)$}
\label{section.M.representation_theory}

\label{subsection.M.simple_modules}

\begin{theorem}
    The monoid $\biheckemonoid=\biheckemonoid(W)$ admits $|W|$
    non-isomorphic simple modules $(S_w)_{w\in W}$ (resp.  projective
    indecomposable modules $(P_w)_{w\in W}$).

    The simple module $S_w$ is isomorphic to the top simple module
    $S^{(w)}_{\{\}}$ of the translation module $T_w$. In general, the
    simple quotient module $S^{(w)}_{J}$ of $T_w$ is isomorphic to
    $S_{\lcoset{w}{J}}$ of $\biheckemonoid$.
\end{theorem}

For example, the simple module $S_{4312}$ is of dimension $3$, with
basis $\{4312,4132,1432\}$ (see Figure~\ref{figure.4312}). The other
simple modules $S_{3412}$, $S_{4123}$, and $S_{1234}$ are of dimension
$5$, $3$, and $1$.

\begin{corollary}
  \label{corollary.translation_module}
  The translation module $T_w$ is an indecomposable
  $\biheckemonoid$-module, quotient of the projective module $P_w$
  of $\biheckemonoid$.
\end{corollary}

$\biheckemonoid_1$ is a submonoid of $\biheckemonoid$. Therefore any
$\biheckemonoid$-module $X$ is a $\biheckemonoid_1$-module, and its
$\biheckemonoid_1$-character $[X]_{\biheckemonoid_1}$ depends only on
its $\biheckemonoid$-character $[X]_{\biheckemonoid}$. This defines a
$\Z$-linear map $[X]_{\biheckemonoid}\mapsto [X]_{\biheckemonoid_1}$.
Let $(S_w)_{w\in W}$ and $(S^1_w)_{w\in W}$ be complete families of
simple modules representatives for $M$ and $\biheckemonoid_1$,
respectively. The matrix of  $[X]_{\biheckemonoid}\mapsto
[X]_{\biheckemonoid_1}$ is called the \emph{decomposition matrix} of
$\biheckemonoid$ over $\biheckemonoid_1$; its coefficient $(u,v)$ is
the multiplicity of $S^1_u$ as a composition factor of $S_v$ viewed as
an $\biheckemonoid_1$-module.
\begin{theorem}
  The decomposition matrix of $M$ over $\biheckemonoid_1$ is upper uni-triangular
  for right order, with $0,1$ entries.
\end{theorem}

\bibliographystyle{alpha}

\bibliography{main}

\end{document}